\documentclass{article}

\usepackage{arxiv}

\usepackage[T1]{fontenc}
\usepackage[utf8]{inputenc}
\usepackage[english]{babel}
\usepackage{microtype}

\usepackage{lmodern}
\usepackage{amsmath}
\usepackage{amssymb}
\usepackage{amsfonts}       % blackboard math symbols
\usepackage{subdepth}
\usepackage{nicefrac}       % compact symbols for 1/2, etc.
\usepackage{dsfont}
\usepackage[longnamesfirst]{natbib}
\bibpunct{(}{)}{;}{a}{,}{,}

\usepackage{url}            % simple URL typesetting
\newcommand*{\doi}[1]{DOI: \href{http://dx.doi.org/\detokenize{#1}}{\detokenize{#1}}}
\usepackage[svgnames]{xcolor}
\usepackage{nameref}
\usepackage{hyperref}
\hypersetup{
  breaklinks,
  colorlinks,
  citecolor={DarkBlue},
  urlcolor={DarkBlue},
  linkcolor={DarkBlue}
}
\usepackage[nameinlink]{cleveref}

\usepackage[textwidth=\marginparwidth,
 backgroundcolor=lightgray,
 linecolor=gray,
 textsize=scriptsize,
 colorinlistoftodos]{todonotes}
\reversemarginpar

\usepackage{graphicx}
\usepackage[ruled,onelanguage]{algorithm2e}
\usepackage{algorithmic}
\usepackage{booktabs}       % professional-quality tables

\newcommand{\pmat}[1]{\begin{pmatrix}#1\end{pmatrix}}

\newcommand{\av}{\ensuremath{\mathbf{a}}}
\newcommand{\bv}{\ensuremath{\mathbf{b}}}

\newcommand{\dv}{\ensuremath{\mathbf{d}}}

\newcommand{\gv}{\ensuremath{\mathbf{g}}}

\newcommand{\kv}{\ensuremath{\mathbf{k}}}

\newcommand{\pv}{\ensuremath{\mathbf{p}}}
\newcommand{\qv}{\ensuremath{\mathbf{q}}}
\newcommand{\rv}{\ensuremath{\mathbf{r}}}
\newcommand{\sv}{\ensuremath{\mathbf{s}}}

\newcommand{\uv}{\ensuremath{\mathbf{u}}}
\newcommand{\vv}{\ensuremath{\mathbf{v}}}
\newcommand{\wv}{\ensuremath{\mathbf{w}}}
\newcommand{\xv}{\ensuremath{\mathbf{x}}}
\newcommand{\yv}{\ensuremath{\mathbf{y}}}
\newcommand{\zv}{\ensuremath{\mathbf{z}}}

\newcommand{\Av}{\ensuremath{\mathbf{A}}}
\newcommand{\Bv}{\ensuremath{\mathbf{B}}}
\newcommand{\Cv}{\ensuremath{\mathbf{C}}}
\newcommand{\Dv}{\ensuremath{\mathbf{D}}}

\newcommand{\Fv}{\ensuremath{\mathbf{F}}}

\newcommand{\Iv}{\ensuremath{\mathbf{I}}}

\newcommand{\Kv}{\ensuremath{\mathbf{K}}}
\newcommand{\Lv}{\ensuremath{\mathbf{L}}}
\newcommand{\Mv}{\ensuremath{\mathbf{M}}}
\newcommand{\Nv}{\ensuremath{\mathbf{N}}}

\newcommand{\Pv}{\ensuremath{\mathbf{P}}}

\newcommand{\Sv}{\ensuremath{\mathbf{S}}}
\newcommand{\Tv}{\ensuremath{\mathbf{T}}}

\newcommand{\Vv}{\ensuremath{\mathbf{V}}}
\newcommand{\Wv}{\ensuremath{\mathbf{W}}}

\newcommand{\Yv}{\ensuremath{\mathbf{Y}}}

\newcommand{\T}{^\mathrm{T}} % Transposed
\newcommand{\mT}{^{-\mathrm{T}}} % Transposed Inverse

\newcommand{\IR}{\mathds{R}} % Reals
\newcommand{\Us}{\mathcal{U}} % Scaled space
\newcommand{\Xs}{\mathcal{X}} % Original space
\newcommand{\C}{^\mathrm{C}}
\newcommand{\xC}{\xv\C} % Cauchy point
\newcommand{\LBFGSB}{\texttt{lbfgsb.m}}
\newcommand{\scaledLBFGSB}{\texttt{scaled-lbfgsb.m}}
\newcommand{\bcflash}{\texttt{Bcflash}}
\newcommand{\cflash}{\texttt{Cflash}}
\newcommand{\scaledbcflash}{\texttt{scaled-Bcflash}}

\title{Scaled Projected-Directions Methods with Application to Transmission Tomography}

\author{
  G.~Mestdagh\\
  Institute of Biomedical Engineering\\
  Polytechnique Montréal\\
  Montreal, QC, Canada \\
  \texttt{guillaume.mestdagh@polymtl.ca}
  \And
  Y.~Goussard\\
  Institute of Biomedical Engineering and Department of Electrical Engineering\\
  Polytechnique Montréal\\
  Montreal, QC, Canada \\
  \texttt{yves.goussard@polymtl.ca}
  \And
  D.~Orban\\
  GERAD and Department of Mathematics and Industrial Engineering\\
  Polytechnique Montréal\\
  Montreal, QC, Canada \\
  \texttt{dominique.orban@polymtl.ca}
}

\begin{document}
\maketitle

\begin{abstract}
Statistical image reconstruction in X-Ray computed tomography yields large-scale regularized linear least-squares problems with nonnegativity bounds, where the memory footprint of the operator is a concern.
Discretizing images in cylindrical coordinates results in significant memory savings, and allows parallel operator-vector products
%\todo{parallel operator-vector products?}
without on-the-fly computation of the operator, without necessarily decreasing image quality.
However, it deteriorates the conditioning of the operator.
We improve the Hessian conditioning by way of a block-circulant scaling operator
% As a consequence, nondiagonal scaling is required to improve the Hessian conditioning.
and we propose a strategy to handle nondiagonal scaling in the context of projected-directions methods for bound-constrained problems.
We describe our implementation of the scaling strategy using two algorithms: TRON, a trust-region method with exact second derivatives, and L-BFGS-B, a linesearch method with a limited-memory quasi-Newton Hessian approximation.
We compare our approach with one where a change of variable is made in the problem.
On two reconstruction problems, our approach converges faster than the change of variable approach, and achieves much tighter accuracy in terms of optimality residual than a first-order method.

\end{abstract}

% keywords can be removed

%%%%%%%%%%%%%%%%%%%%%%%%%%%%%
% BODY
%%%%%%%%%%%%%%%%%%%%%%%%%%%%%
\newpage
%%%%%%%%%%%%%%%%%%%%
\section{Introduction}

% General problem
We consider the bound-constrained problem
\begin{equation}\label{eq:or-pb}
    \min\ f(\xv) \quad \mbox{s.t.}\ \xv \geqslant 0,
\end{equation}
where $f:\IR^n \to \IR$ is convex and $\mathcal{C}^2$.
We assume that $\nabla^2 f$ cannot be stored explicitly or in factorized form.
%
% Properties of the problem
We are particularly interested in the case where \eqref{eq:or-pb} is large and badly scaled.
Our main motivation is to solve efficiently statistical image reconstruction problems arising from X-Ray Computed Tomography (CT) \citep{Herman09}.
Whereas cartesian coordinates are typical, discretizing such problems in cylindrical coordinates yields large savings in storage, but results in badly scaled problems and, without proper scaling, off-the-shelf solvers usually fail.

% Proposed solution
In this paper we employ a scaling strategy that exploits the structure of \(f\) combined with the trust-region projected Newton method of \cite{Lin99} and with the line search limited-memory BFGS for bound-constrained problems of \cite{Byrd95} to maintain satisfaction of the bound constraints at all times.

\subsection*{Motivation and previous work}

% Introduction
Our main interest resides in statistical image reconstruction problems arising from X-Ray Computed Tomography (CT) \citep{Herman09}.
Compared to analytical methods such of the filtered backprojection family \citep{Feldkamp84}, statistical reconstruction results in less noisy images but is more computationally expensive \citep{Fessler00}.

% Reconstruction problem statement
\citet{Sauer93} show that an image $\xv$ can be estimated from the measurements $\bv$ by solving
\begin{equation}\label{eq:recmodel}
	\min\ \tfrac{1}{2} \|\Av \xv - \bv \|^2_\Vv + \lambda R(\xv)
	\quad\mbox{s.t.}\ \xv \geqslant 0,
\end{equation}
where $\Av$ is a large sparse matrix, $\lambda > 0$ is a regularization parameter, \(R: \IR^n \to \IR\) is a convex regularization function and $\Vv$ is a diagonal weight matrix with $V_{ii} = \exp(-b_i)$.
In \eqref{eq:recmodel}, the objective is composed of a least-squares data-fitting term, and a regularization term to discourage large differences between adjacent pixels.
Here, we focus on situations where one wishes to solve problem~\eqref{eq:recmodel} precisely (with a low tolerance level), which corresponds to instances where the information present in the data should be fully exploited, as, e.g., in low-dose CT reconstruction or for imaging 
thin structures in micro-CT \citep{Hamelin10a}.

%% Discretization in cylindrical coordinates
While reconstructed images are usually discretized using a cartesian grid of voxels \citep[equation (6)]{Fessler00}, we use cylindrical coordinates.
This discretization is in adequation with the circular geometry of the data acquisition process and results in block-circulant $\Av$, which has far lower storage requirements than the operator resulting from cartesian coordinates.
Moreover, it is well suited for parallelization and does not require computing the projection matrix on the fly \citep{Goussard13}.
\citet{Thibaudeau13} present details about projection matrix storage and cylindrical-to-cartesian conversion to take advantage of the projection operator structure without loss of resolution.

%% State-of-the art in image reconstruction
In statistical image reconstruction, first-order methods are often preferred due to their low storage and computational demands, and their simplicity \citep{Fessler00}.
Usual reconstruction methods include the expectation-maximization algorithm \citep{Lange95,Ahn06}, coordinate-descent methods \citep{Sauer93, Noo16} and gradient-based methods \citep{Erdogan99,Kim14}.
This last category is probably the most studied in image reconstruction.
Such methods are often improved by using ordered subsets \citep{Erdogan99a,Hudson94} or Nesterov momentum \citep{Nesterov83,Kim14,Xu16,Choi10,Jensen12}.
Recently, problem-splitting and proximal approaches have been proposed in the context of sparse reconstruction \citep{Sidky13,Nien15,Xu16}.

Among the most popular first-order algorithms is the spectral projected gradient method (SPG) of \citet{Birgin02}. 
This variant of the projected gradient method involves a Barzilai-Borwein steplength \citep{Barzilai88} and a linesearch to ensure convergence.
Due to its simplicity and efficiency, SPG is used for inverse problems in engineering applications \citep[section~4]{Birgin14}, including image deblurring \citep{Bonettini08} and radiology \citep{Kolehmainen06,Kolehmainen07}.

For imaging applications, diagonal scaling is usually sufficient for first-order methods to perform well \citep{Pock11,Bonettini08}.
In CT reconstruction, \citet{Erdogan99a} propose a diagonal scaling strategy referred to as the \textit{separable quadratic surrogate method}.
This strategy is used in several reconstruction methods \citep{Kim14, Nien15, Zheng18}.
\citet{Piccolomini18} use the diagonally scaled gradient method of \citet{Bonettini08} to reconstruct CT images in the context of incomplete data.
However, in cylindrical coordinates, widely different voxel sizes make \(\Av\) badly scaled, and diagonal scaling is no longer appropriate.

%% Block-circulant structures and scaling matrix design
To improve the conditioning of \eqref{eq:recmodel}, we follow \cite{Golkar13} and use a block-circulant scaling operator $\Cv$ that exploits the block-circulant property of $\Av$ and of the finite-difference matrices that appear in the regularization term $R(\xv)$.
The scaling operator can be written
\begin{equation*}
    \Cv = \Fv^{\star} \Tv \Fv,
\end{equation*}
where $\Tv$ is diagonal, $\Fv$ is a discrete block Fourier transform, and a star indicates the conjugate transpose.
Thus, $\Cv$ and its inverse can be applied to a vector at the cost of a fast Fourier transform, namely in $\mathrm{O}(s_b n_b \log n_b)$ operations, where $s_b$ is the size of a square block, and $n_b$ is the number of blocks.
To obtain the coefficients of $\Tv$, we block-diagonalize the Hessian $\nabla f(\xv)$ (or its approximation when necessary) and use the inverse diagonal coefficients of the resulting block-diagonal matrix. 
In other words
\begin{equation*}
       \Tv = \mathrm{diag}\left(\Fv \nabla f(\xv) \Fv^\star\right)^{-1}
\end{equation*}

Note that, though matrix-vector products with $\Cv$ are computed using a block Fourier transform, it is not the case for $\Av$.
The matrix $\Fv \Av \Fv^\star$ is block-diagonal, but its diagonal blocks are dense, and storing all of them is too costly and would defeat the purpose of using cylindrical coordinates.
For this reason, we store the first line of blocks of $\Av$
using an incremental compressed-columns-sparse (ICCS) scheme.
We compute the diagonal blocks of $\Fv \Av \Fv^\star$ one at a time to compute the diagonal of $\Tv$.
As a consequence, computing products with $\Av^{-1}$ requires using iterative methods.

The change of variable $\xv = \Cv \uv$ transforms~\eqref{eq:or-pb} into the scaled problem
\begin{equation}\label{eq:screcmodel}
	\min \ \tfrac{1}{2} \|\Av \Cv\uv - \bv \|^2_\Vv + \lambda R(\Cv\uv)
	\quad\mbox{s.t.}\ \Cv\uv \geqslant 0.
\end{equation}
If $\Cv$ satisfies $\Cv\T \nabla^2 f(\xv) \Cv \approx \Iv$, the objective Hessian is better conditioned in \eqref{eq:screcmodel} than in \eqref{eq:recmodel}.
However, \eqref{eq:screcmodel} features linear inequality constraints instead of simple bounds.

% Maxime
\citet{McLaughlin17} solves \eqref{eq:screcmodel} with \cflash{}, a variant of TRON \citep{Lin99} adapted to linear inequality constraints.
Each iteration of \cflash{} requires projecting candidate iterates \(\uv\) into the feasible set by solving
\begin{equation}
   \label{eq:proj-into-Cv>=0}
% 	\mathrm{Proj}(\uv) \in \arg
	\min_{\vv}\ \|\vv - \uv\| \quad \mbox{s.t.}\ \Cv \vv \geqslant 0,
\end{equation}
which represents a significant amount of computation.
In \cflash{}, the above projection is computed efficiently by solving the dual problem, which is a bound-constrained linear least-squares problem with operator \(\Cv\), iteratively.
Even though~\eqref{eq:proj-into-Cv>=0} can be solved efficiently thanks to the structure of $\Cv$, it remains substantially more costly than projecting onto simple bounds.

Instead of solving \eqref{eq:screcmodel}, we propose to solve (\ref{eq:recmodel}) with a scaled quasi-Newton and Newton method in order to reproduce the effect of a change of variable without actually performing it.
Thus, we benefit from the conditioning improvement provided by the scaling operator, yet we still only need to perform projections onto the nonnegative orthant via the direct formula
\begin{equation*}
    \mathrm{Proj}(\xv) = \max(\xv, 0),
\end{equation*}
where the $\max$ applies componentwise.
There are other advantages to using scaled directions, including the possibility to choose a new scaling operator at each iteration.
In this, work, though, we use the same scaling for all iterations.
Moreover, the optimality measure is the same whatever the scaling, which makes it possible to compare easily different scaling operators in terms of convergence speed.
Because there is no change of variable, the coordinates of an iterate are the pixels of the image. It is possible to access them at no cost during the optimization process, for instance to visualize the progress of the reconstruction.

% The Italians
\citet{Bonettini08} describe a diagonally-scaled variant of SPG in which both the gradient and projection subproblems are scaled.
\cite{Bonettini13} extend the approach to block-diagonal scaling in the context of image deblurring, a problem where the system operator is a 2D convolution, which is block-circulant with circulant blocks.
Their scaling depends on the active constraints at the current iterate, as in the projected Newton method described by \cite{Bertsekas82}, which allows them to apply a nondiagonal scaling while preserving simple projections.
They also investigate a quasi-Newton method where the Hessian approximation is a truncated spectral decomposition of the system matrix.
In the partially-diagonal scaling approach, \citeauthor{Bonettini13} apply the method of \cite{Landi2008} and solve the linear system inexactly using the conjugate gradient method (CG) of \cite{Hestenes52}.
Both methods are presented as scaled gradient projection methods in which the scaling operator is inspired from the problem Hessian or from a quasi-Newton matrix.

% Us compared to the Italians
Our approach extends that of \citet{Bonettini08} to more complex methods for bound-constrained problems.
Because we apply the scaling to higher-order methods, we restrict our study to the situation where the scaling operator can be applied or inverted easily.
The problem Hessian does not enter this category, because it is expensive to apply and we need to use CG to solve a linear system involving it.
While \citeauthor{Bonettini13} apply a partially-diagonal mask to the problem Hessian, we apply the partially-diagonal mask on the fast scaling operator or its inverse,
in the context of a change of scalar product, which is compatible with any solution method. %independently from the solution method.
We then describe the impact of the scaling on families of general-purpose optimization algorithms.
Specifically, we consider the limited-memory BFGS method for bound-constrained problems, and a trust-region projected Newton method.
The convergence properties of both methods rest on the computation of a Cauchy point, i.e., an approximate minimizer in the negative gradient direction.
Our implementation of L-BFGS-B differs from that of \cite{Zhu97} in that we compute an inexact Cauchy point and restrict the computation of a step to an active face by way of restriction operators.
Our implementation of a projected Newton method follows the design of TRON \citep{Lin99}, in which projected gradient steps are performed to identify an inexact Cauchy point and a candidate active set, followed by a sequence of Newton steps on the active face of the feasible set globalized by a trust-region mechanism.
As opposed to a traditional active-set method, which only adds or removes one bound at a time from the active set estimate, a projected direction approach has been shown to be able to add and remove many bound constraints from the active set at a time and to be particularly appropriate for large-scale problems.
We illustrate the performance of both methods on synthetic images, and we compare it to that of the projected spectral gradient method, a classic method in imaging.

%% Notation
\paragraph{Notation}
Lowercase and uppercase bold Latin letters represent vectors and matrices, respectively. 
Light face Latin letters represent integers, such as iteration counters, and functions.
In addition, the \(i\)-th component of \(\xv\) is denoted \(x_i\) and the \((i,j)\)-th element of \(\Av\) is \(A_{i\!j}\).
Light face Greek letters represent scalars.
The Euclidean scalar product on $\IR^n$ is denoted $\langle \cdot, \cdot\rangle$.
The $i$-th partial derivative of function $f$ at $x$ is denoted \(\partial_i f(x)\).

%%%%%%%%%%%%%%%%%%%%
\section{A scaling strategy for bound-constrained problems}

In this section we describe the effect of scaling on procedures that are common between the two methods we present.
After a short presentation of the scaling strategy, we present two procedures we use in L-BFGS-B and TRON to compute a Cauchy point and a descent step in the active face respectively.
We then apply the scaling on the limited-memory quasi-Newton matrix that appears in L-BFGS-B.
We integrate these procedures into the chosen algorithms in \Cref{S3}.

\label[section]{S2}

    \subsection{Overview of the strategy}
    	Our scaling strategy consists in using a metric in which the problem is well scaled and the projections can be computed with a direct formula.

A linear change of variables is equivalent to changing the scalar product in the original space.
Indeed, if $\xv = \Cv \uv$ and $\zv = \Cv \wv$,
\begin{equation*}
    \langle \uv, \wv  \rangle = \langle \xv, \Pv^{-1} \zv\rangle,
    \qquad
    \Pv := \Cv \Cv\T.
\end{equation*}
This equivalence makes it possible to import geometric elements from the scaled space into the original space.
From now on, we use $\Xs$ to denote the original space and $\Us$ for the scaled space.
Every $\xv\in \Xs$ corresponds to a $\uv\in\Us$ such that $\xv = \Cv \uv$.
We denote $\bar f$ the objective function of \eqref{eq:or-pb} in the scaled space, i.e., for $\uv \in \Us$,
\begin{equation*}
    \bar f(\uv) := f(\Cv \uv),
    \qquad
    \nabla \bar f(\uv) = \Cv\T \nabla f(\Cv \uv),
    \qquad
    \nabla^2 \bar f = \Cv\T \nabla^2 f(\Cv \uv) \Cv.
\end{equation*}

The first element we transform is the gradient direction.
Indeed, due to the choice of scaling, the gradient direction in $\Us$ is expected to be a more promising descent direction than the gradient direction in $\Xs$.
If $\uv \in \Us$, $\xv = \Cv\uv$, and $\alpha > 0$,
\begin{equation*}
    \xv(\alpha) :=
    \Cv (\uv - \alpha \nabla \bar f (\uv) ) 
    = \xv - \alpha \Cv \nabla \bar f (\uv) 
    = \xv - \alpha \Pv \nabla f (\xv). 
\end{equation*}
In other words, a step along the negative gradient in $\Us$ is equivalent to a step in the direction 
\begin{equation}\label{eq:scaledgradient}
	\qv = -\Pv\nabla f (\xv)
\end{equation} 
in $\Xs$.
We use~\eqref{eq:scaledgradient} instead of $-\nabla f(\xv)$ in the hope that the scaled search direction better captures natural problem curvature.

    \subsection{Projected directions} 
    	\label{proj}
In the context of projected methods, it is often necessary to work on a face of the feasible set, or to take projected gradient steps.

\subsubsection{Projected gradient steps}
\label{sec:proj-methods}

A standard projected gradient step from $\xv$ can be described by \(\mathrm{Proj}\left(\xv - \alpha \nabla f(\xv) \right) - \xv\) where $\alpha > 0$.
A scaled projected gradient step has the form \(\xv(\alpha) - \xv\), where
\begin{equation*}
	\xv(\alpha) = \mathrm{Proj}\left(\xv + \alpha \dv \right),
\end{equation*}
and where \(\dv\) is a linear transformation of $\nabla f(\xv)$.
The direction $\dv$ must be a descent direction in the sense that there exists $\bar\alpha > 0$ such that $f(\xv(\alpha)) < f(\xv)$ for all $\alpha \in \left(0, \bar\alpha\right]$.

\citet[section 2]{Bertsekas82} explains that~\eqref{eq:scaledgradient} might not be such a descent direction.
We define the \textit{binding} constraints at $\xv$ as those with indices in
\begin{equation}\label{eq:bindingset}
   I_+(\xv) = \{i\ |\ x_i = 0\mbox{ and }\partial_i f(\xv) > 0\},
\end{equation}
where $\partial_i f(\xv)$ is the $i$-th component of $\nabla f(\xv)$.
We introduce the subspace
\begin{equation*}
	F = \{\av \in \IR^n \mid a_i = 0 \text{ for all } i \in  I_+(\xv)\},
\end{equation*}
and the set $F_+ = F \cap \IR_+^n$, called the face of the feasible set exposed by \(-\nabla f(\xv)\).
\citet[Proposition~1]{Bertsekas82} establishes that
\begin{equation}\label{eq:desc}
    \dv = -\bar \Pv \nabla f(\xv),
\end{equation}
where the matrix $\bar \Pv$ is defined by
\begin{equation}\label{eq:pbar}
        \bar\Pv_{ij} = \left\{\begin{array}{ll}
        \Pv_{ij} & \mbox{if } i, j \notin I_+(\xv) \\
        0 & \mbox{otherwise,} 
    \end{array}\right.
\end{equation}
is a descent direction in the sense defined above.

To compute (\ref{eq:desc}), we decompose $\IR^n = F \oplus G$ where $G = F^\perp$, and we write
\begin{equation*}
	\nabla f(\xv) = \pmat{\gv_F \\ \gv_G}\qquad
	\Pv = \pmat{\Pv_{F\!F} & \Pv_{F\!G} \\ \Pv_{G\!F} & \Pv_{G\!G}},
\end{equation*}
where $\gv_F = (\partial_i f(\xv))_{i\notin I_+(\xv)}$ and $\gv_G = (\partial_i f(\xv))_{i\in I_+(\xv)}$.

The gradient is first projected onto the subset $F$.
Then the scaling is made on the projected gradient $\gv_F$ by applying the principal submatrix $\Pv_{FF}$.
This submatrix is obtained by keeping only the rows and columns whose indices are not in $I_+(\xv)$.
Finally,
\begin{equation*}
	\dv = -\pmat{\Pv_{F\!F} & 0 \\ 0 & 0}\pmat{\gv_F \\ \gv_G} = \pmat{-\Pv_{F\!F\,}\gv_F \\ 0}.
\end{equation*}

\subsubsection{Conjugate gradient in a face of the feasible set}\label{cg}

The same procedure can be used to modify conjugate gradient directions inside a face of the feasible set.
Consider the quadratic problem
\begin{equation}\label{eq:quadinF}
	\min_{x\in F} \tfrac{1}{2} \xv\T \Bv \xv + \xv\T \gv,
\end{equation}

To solve~(\ref{eq:quadinF}), the Conjugate Gradient method of \citet{Hestenes52} is applied to the equivalent reduced problem
\begin{equation*}
	\min_{\bar{x}\in \IR^{\dim F}} \tfrac{1}{2} \bar{\xv}\T \Bv_{FF} \bar{\xv} + \bar{\xv}\T \bar\gv.
\end{equation*}

The directions $\bar{\pv}_1, \bar{\pv}_2, \dots$ generated by the procedure are conjugate with respect to the principal submatrix $\Bv_{FF}$ of \(\Bv\) \citep{Hestenes52}.
In particular, at the $k$-th iteration, the next direction is defined as
\begin{equation*}
	\bar\pv_{k+1} = \bar\rv_{k+1} + \beta_k \bar\pv_k,
\end{equation*}
where $\bar{\rv}_{k+1} = \bar\gv - \Bv_{FF} \bar\xv_{k+1}$ is the residual and $\beta_k$ is chosen so that the conjugacy condition $\bar{\pv}_{k+1}\T \Bv_{FF} \bar{\pv}_{k} = 0$ is verified.

To improve the convergence of CG, we use a scaled residual to generate the new direction. The direction update formula becomes
\begin{equation*}
	\bar\pv_{k+1} = \Pv_{FF} \bar\rv_{k+1} + \beta'_k \bar\pv_k,
\end{equation*}
where $\beta'_k$ is chosen to respect the conjugacy condition.
Note that applying the scaling in the case of CG, is equivalent to preconditioning CG with $\Pv_{FF}$.

    \subsection{Limited memory quasi-Newton matrices}
    	\label{qn}

In a quasi-Newton method, the objective function is approximated about the current iterate $\xv_k$ by the quadratic model
\begin{equation}\label{eq:or-qnmodel}
	m_k(\xv_k  + \zv) = f_k + \gv_k\T \zv + \tfrac{1}{2}\zv\T \Bv_k \zv,
\end{equation}
where $f_k$ and $\gv_k$ are the objective value and gradient at $\xv_k$ respectively,
and $\Bv_k = \Bv_k^T$ is an approximation of $\nabla^2 f(\xv_k)$.
Secant methods are a special case in which \(\Bv_k\) must satisfy the secant equation
\begin{equation*}
    \Bv_k \sv_{k-1} = \yv_{k-1},
\end{equation*}
where \(\sv_{k-1} := \xv_k - \xv_{k-1}\) and \(\yv_{k-1} := \nabla f(\xv_k) - \nabla f(\xv_{k-1})\).
Secant methods typically define \(\Bv_k\) as rank-one or rank-two update of \(\Bv_{k-1}\) involving \(\sv_{k-1}\) and \(\yv_{k-1}\).
This has the disadvantage that \(\Bv_k\) is almost always dense even though \(\nabla^2 f(\xv_k)\) might be sparse.
Therefore, the entire matrix \(\Bv_k\) must be stored, which is unrealistic in large-scale applications.
However, at least conceptually, a product between \(\Bv_k\) and a vector could be computed without storing \(\Bv_k\) if the initial approximation \(\Bv_0\) along with all the pairs \(\{\sv_i, \yv_i\}_{0 \leq i \leq k-1}\) are stored instead.

In a limited-memory context, we store an initial matrix $\Bv_0$ along with the $m$ most recent pairs $\{\sv_i, \yv_i\}_{k-m \leqslant i \leqslant k-1}$, where $m$ is the memory.
The procedure uses the information from $\Bv_0$ and from the \(m\) pairs to update and compute a product with $\Bv_k$.
Even though \(\Bv_k\) would still be dense if it were materialized, it it only represented implicitly.

Because the memory $m$ is often small compared to the problem dimension, quasi-Newton updates can only contribute a limited amount of information to $\Bv_k$.
For this reason, the choice of $\Bv_0$ is critical to obtain a good approximation of the objective Hessian.
In particular, when the Hessian is ill conditioned, choosing $\Bv_0$ as a multiple of the identity might lead to poor performance.

The best-known limited-memory quasi-Newton method is probably the limited-memory BFGS method \citep{Nocedal80}, which additionally ensures that \(\Bv_k\) is positive definite provided that \(\Bv_{k-1}\) is positive definite and \(\sv_{k-1}^T \yv_{k-1} > 0\).
Several procedures exist to compute a product between \(\Bv_k\) or its inverse and a vector, including the two-loop recursion \citep{Nocedal80}, and a variant based on compact storage \citep{Byrd94a}.
We present the second one because it was designed to handle bound constraints.

The pairs $\{\sv_i, \yv_i\}_{k-m \leqslant i \leqslant k-1}$ are stored in two matrices 
\begin{equation*}
	\Sv = \pmat{\sv_{k-m} & \cdots & \sv_{k-1}}\quad\mbox{and}\quad
	\Yv = \pmat{\yv_{k-m} & \cdots & \yv_{k-1}},
\end{equation*}
and we define
\begin{equation*}
	\Dv = \begin{pmatrix}
		\sv\T_{k-m} \yv_{k-m} & & \\
		 & \ddots & \\
		 & & \sv\T_{k-1} \yv_{k-1}
		\end{pmatrix}
		\quad\text{and}\quad
		\Lv = \begin{pmatrix}
		0 & \cdots & \cdots & 0 \\
		\sv_{k-m+1}\T \yv_{k-m} & \ddots & & \vdots\\
		\vdots & \ddots & \ddots & \vdots \\
		\sv_{k-1}\T \yv_{k-m} & \cdots & \sv_{k-1}\T \yv_{k-2} & 0
	\end{pmatrix}.
\end{equation*}
In the compact formula, \(\Bv_k\) is stored implicitly as \(\Bv_0\),
\begin{equation}\label{eq:WandM}
	\Wv = \begin{pmatrix}
 	     \Yv & \Bv_0 \Sv
	     \end{pmatrix}
	\quad \text{and} \quad  
	\Mv = \begin{pmatrix}
      	-\Dv & \Lv\T           \\
      	\Lv  & \Sv\T \Bv_0 \Sv
      	\end{pmatrix}^{-1},
\end{equation}
such that
\begin{equation}\label{eq:compactformula}
	\Bv_k = \Bv_0 - \Wv \Mv \Wv\T.
\end{equation}
where $\Bv_0$ is positive definite.

In most implementations, $\Bv_0$ is chosen as
\begin{equation}\label{eq:theta}
	\Bv_0 = \theta \Iv\quad\mbox{with}\quad	
	\theta = \frac{\yv_{k-1}\T \yv_{k-1}}{\yv_{k-1}\T \sv_{k-1}},
\end{equation}
where $\theta$ is a scaling parameter \citep{Byrd95}.
A diagonal $\Bv_0$ leads to very efficient operations with the L-BFGS compact formula.
However, it might be inappropriate for approximating ill-conditioned Hessians.

We choose $\Bv_0$ so that the L-BFGS operator in $\Xs$ reproduces the behavior of a standard L-BFGS operator with initial approximation~\eqref{eq:theta} in $\Us$.

Assume that, in \(\Us\), $\bar{f}$ is approximated about the current scaled iterate $\uv_k$ by the quadratic model
\begin{equation}\label{eq:sc-qnmodel}
	m'_k(\uv_k + \wv) = \bar{f}(\uv_k) + \nabla \bar{f}(\uv_k)\T \wv
              + \tfrac{1}{2} \wv\T \Bv_k' \wv,
\end{equation}
where $\Bv_k'$ is a L-BFGS operator with initial approximation~\eqref{eq:theta}.
The pairs $\{\bar\sv_i, \bar\yv_i\}$ in \(\Us\) are related to the pairs $\{\sv_i, \yv_i\}$ in \(\Xs\) via
\begin{equation}\label{eq:relpair}
	\bar\sv_i = \Cv^{-1} \sv_i\qquad
	\bar\yv_i = \Cv\T \yv_i,\qquad
	i = k-m, \dots, k-1.
\end{equation}
We replace $\{\sv_i, \yv_i\}$ with $\{\bar\sv_i, \bar\yv_i\}$ in \eqref{eq:compactformula} and \eqref{eq:theta}, and obtain
\begin{equation}\label{eq:scaledcompact}
	\Bv_k' = \bar\theta \Iv - 
	\begin{pmatrix}
		\Cv\T \Yv && \bar\theta \Cv^{-1} \Sv
	\end{pmatrix}
	\begin{pmatrix}
		-\Dv && \Lv\T                            \\
		\Lv  && \bar\theta\Sv\T \Pv^{-1} \Sv
	\end{pmatrix}^{-1}
	\begin{pmatrix}
		\Yv\T \Cv          \\
		\bar\theta  \Sv \Cv\mT
	\end{pmatrix},
\end{equation}
where
\begin{equation}\label{eq:omega}
	\bar\theta = \frac{\bar\yv_{k-1}\T \bar\yv_{k-1}}{\bar\sv_{k-1}\T \bar\yv_{k-1}}
	       = \frac{\yv_{k-1}\T \Pv \yv_{k-1}}{\yv_{k-1}\T \sv_{k-1}}.
\end{equation}

We use (\ref{eq:or-qnmodel}) to approximate $f$ in \(\Xs\).
A comparison between (\ref{eq:or-qnmodel}) and (\ref{eq:sc-qnmodel}) yields
\begin{equation*}
	\Bv_k' = \Cv\T \Bv_k \Cv.
\end{equation*}
Finally $\Bv$ is a L-BFGS operator with initial approximation
\begin{equation}\label{eq:B0}
	\Bv_0 = \bar\theta \Cv\mT \Cv^{-1} = \bar\theta \Pv^{-1}
\end{equation}

Apart from the storage of $\Pv$, this modification does not require more storage in the compact L-BFGS formula.
Instead of storing $\Sv$ and $\Sv\T \Sv$, we store $\Pv^{-1} \Sv$ and $\Sv\T \Pv^{-1} \Sv$, while $\Lv$ and $\Dv$ remain unchanged.
The L-BFGS update only requires one product with $\Pv^{-1}$ to compute $\Pv^{-1} \sv_k$ and one scalar product defined by $\Pv$ to compute $\bar\theta$.

%%%%%%%%%%%%%%%%%%%%
\section{Modified algorithms}\label[section]{S3}

In this section, we present the salient elements of the L-BFGS and TRON algorithms, and of our implementations.
Then we describe the modification we brought to apply the scaling strategy.

    \subsection{The L-BFGS-B algorithm}
    	The L-BFGS-B algorithm of \citet{Byrd95} is a popular quasi-Newton method for bound-constrained problems.
Its standard implementation exploits the compact representation of limited-memory quasi-Newton operators, a diagonal \(\Bv_0\), and the Sherman-Morrison-Woodbury formula to solve linear systems whose coefficient is a principal submatrix corresponding to inactive indices.
In our application, \(\Bv_0\) is nondiagonal and its principal submatrices are not structured, so the Sherman-Morrison-Woodbury approach would be inefficient.

\subsubsection{Presentation of the algorithm}

At the beginning of an iteration, we compute the Cauchy point \(\xv_k^\mathrm{C}\) as the \emph{exact} first local minimizer of (\ref{eq:or-qnmodel}) along the piecewise affine path
\begin{equation}\label{eq:projgradpath}
	t \mapsto \mathrm{Proj}(\xv_k - t\,\gv_k),
\end{equation}
where $\gv_k$ is the objective gradient at the current iterate $\xv_k$.
The Cauchy point is obtained by successively examining the quadratic model~\eqref{eq:or-qnmodel} on each segment of (\ref{eq:projgradpath}). 
On a segment between two breakpoints, the model is a second-order polynomial function of the nonnegative parameter $t$. 
If the polynomial is nonincreasing on the segment, then the procedure moves to the next segment. 
Otherwise a minimizer is computed on the current segment and returned as the Cauchy point~\citep{Byrd95}.

For clarity, we now drop the iteration index \(k\).
The Cauchy point is used to identify the set of active constraints
\begin{equation}
	\mathcal{A} = \left\{i\ |\ x_i^\mathrm{C} = 0\mbox{ and }g_i > 0\right\}.
\end{equation}
This characterization is a consequence of the procedure presented by \citet[Section 4, Algorithm CP]{Byrd95}.
Their algorithm returns a set of free indices $\mathcal{F}$, which is the complementary of $\mathcal{A}$.
The active constraints are associated to breakpoints that are between $\xv$ and $\xv^\mathrm{C}$ along the projected gradient path $\eqref{eq:projgradpath}$.
Note that $\mathcal{A}$ is not a binding set as defined in~\eqref{eq:bindingset} because its definition involves $\gv = \nabla f(\xv)$, i.e., the objective gradient is not evaluated at $\xv^\mathrm{C}$.
Using $\mathcal{A}$ we define the subspace
%\todo{DO: n'est-il pas suffisant de dire $x_i = 0$ pour $i \in A$ ?}
\begin{equation}\label{eq:activeface}
	F = \left\{\xv\in\IR^n\ |\ \forall\, i \in \mathcal{A}\quad x_i = x_i^\mathrm{C} = 0\right\},
\end{equation}
and the active face $F_+ = F \cap \IR_+^n$.
In order to compute a minimizer of~\eqref{eq:or-qnmodel} over $F$ and obtain a descent direction, we set $\Bv = \Bv_k$ and $\gv = \gv_k$ in~\eqref{eq:quadinF} and solve by way of the Sherman-Morrison-Woodbury formula.
If we partition
\begin{equation*}
	\Wv = \pmat{\Wv_F \\ \Wv_G},
\end{equation*} 
where $\Wv$ is defined in (\ref{eq:WandM}), the inverse of a principal submatrix of \(\Bv\) is
\begin{equation}\label{eq:smw}
	\Bv_{F\!F}^{-1} = \Bv_{0,F\!F}^{-1}
	              - \Bv_{0,F\!F}^{-1} \Wv_F
	              \left(\Mv^{-1} - \Wv_F\T \Bv_{0,F\!F}^{-1} \Wv_F \right)^{-1}
	              \Wv_F\T \Bv_{0,F\!F}^{-1}.
\end{equation}
With $\Bv_0 = \theta \Iv$, we have
\begin{equation*}\label{eq:smw2}
	\Bv_{F\!F}^{-1} = \theta^{-1} I
	              - \theta^{-2} \Wv_F
	              \left(\Mv^{-1} - \theta^{-1} \Wv_F\T\Wv_F \right)^{-1}
	              \Wv_F\T.
\end{equation*}

Finally, we use a strong Wolfe linesearch
% between the current iterate and the solution of (\ref{eq:quadinF})
to determine the next iterate \(\xv_{k+1} = \xv_k + \alpha \dv\), where \(\alpha > 0\) and \(\dv\) is the search direction from~\eqref{eq:quadinF}.
The procedure, proposed by \citet{More94}, returns, if possible, a steplength \(\alpha\) such that \(x_{k+1}\) satisfies the bound constraints and the conditions
\begin{equation}\label{eq:strwlf}
    f(\xv_{k+1}) \leqslant f(\xv_k) + \mu \alpha\, \gv_k\T \dv
    \quad\mbox{and}\quad
    |\gv_{k+1}\T \dv| \leqslant \eta |\gv_k\T \dv|,
\end{equation}
where $0 < \mu < \eta < 1$ are parameters.
\Cref{alg:lbfgsb} shows an overview of the L-BFGS-B method.

\begin{algorithm}[h]\label[algorithm]{alg:lbfgsb}
    \caption{Overview of the standard L-BFGS-B algorithm.}
    \KwData{$\xv_0$, parameters}
    \For{$k = 0, 1, 2, \dots$}{
        Compute an exact Cauchy Point $\xv_k^\mathrm{C}$ along the projected path
        $t \mapsto \mathrm{Proj}\left(\xv_k - t\,\nabla f(\xv_k)\right)$ \\
        Identify the active face $F_+$\\
        Find a minimizer of the model (\ref{eq:or-qnmodel}) over the affine subspace
        (\ref{eq:activeface}) by the Sherman-Morrison-Woodbury formula \\
        Perform a strong Wolfe linesearch to find the next iterate $\xv_{k+1}$ satisfying \eqref{eq:strwlf}\\
        Update the L-BFGS operator.
    }
\end{algorithm}

\subsubsection{Implementation in Matlab}

Our MATLAB implementation of L-BFGS-B\footnote{Available online at \url{https://github.com/optimizers/NLPLab}}, \LBFGSB{}, solves~\eqref{eq:quadinF} with CG instead of the Sherman-Morrison-Woodbury formula.
Indeed, though \LBFGSB{} uses \eqref{eq:theta} as an initial matrix, our implementation should work with \eqref{eq:B0} and \eqref{eq:omega}.
When \(\Bv_0\) is nondiagonal, computing one of the products between \(\Bv_{0,F\!F}^{-1}\) and a vector in \eqref{eq:smw} requires to solve a linear system and might be as expensive as computing a product with \(\Bv_{F\!F}^{-1}\) by CG.

\begin{figure}
    \centering
    \includegraphics[height=.31\columnwidth]{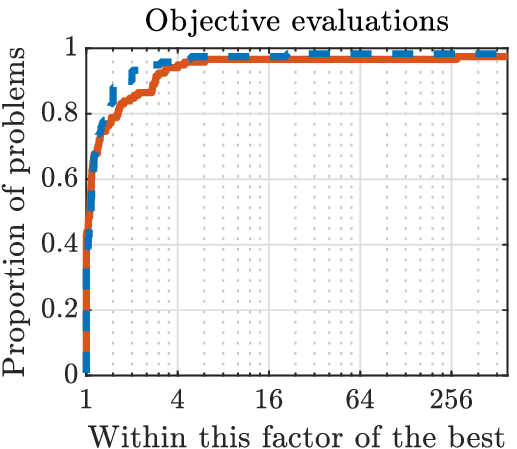}\!
    \includegraphics[height=.31\columnwidth]{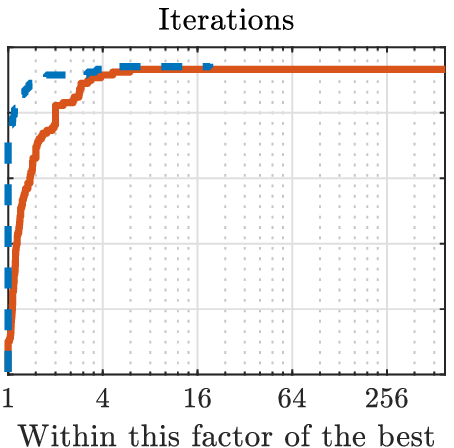}\!
    \includegraphics[height=.31\columnwidth]{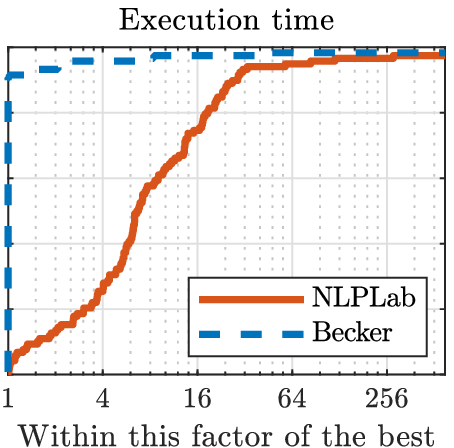}\!
    \caption{Log-scale performance profiles \citep{Dolan02} of our MATLAB implementation versus the C interface. The values compared are the number of objective evaluations~(left), the number of iterations~(middle) and the execution time~(right).}
    \label[figure]{fig:valLBFGSB}
\end{figure}

We validated our implementation against a C translation of the original Fortran implementation \citep{Zhu97,Morales11} provided by Stephen Becker\footnote{Available online at \url{https://github.com/stephenbeckr/L-BFGS-B-C}} on a collection of standard problems.
The C version uses the Sherman-Morrison-Woodbury formula to solve~\eqref{eq:quadinF}.
Our benchmark comprises $128$ bound-constrained problems from the CUTEst library \citep{Gould15}.
The tests ran on a 3GHz Intel Core i7-5960X with 64GB of RAM\label{compspec}.
We report our results in the form of performance profiles \citep{Dolan02} in logarithmic scale in \Cref{fig:valLBFGSB}.
On the performance profile, a point with coordinates $(x, y)$ on a solver curve means that, on a proportion $y$ of the problems, the solver's performance measure was at most $x$ times that of the best solver.

Because we use MATLAB instead of C, \LBFGSB{} is slower than the original version.
However, the results are similar in terms of number of objective evaluations and iterations, even though a slight degradation in performance on standard problems is somewhat expected and matches the observations of \cite{Byrd95}.

\subsubsection{Modifications related to scaling}

In order to obtain a better Hessian approximation, we use (\ref{eq:omega}) and (\ref{eq:B0}) in the L-BFGS operator.
A non-diagonal $\Bv_0$ requires several modifications in the algorithm, as the procedures presented by \citet{Byrd95} owe their efficiency to the diagonal structure of $\Bv_0$.

Finding the Cauchy point requires examining up to $n$ segments defined by the breakpoints along the projected gradient path.
The feasibility of an exact search relies on the absence of any operation with complexity worse than $O(n)$ in the update of the quadratic model derivatives along each segment. 
In particular, the scalar product between a row of $\Bv_0$ and a vector is required for each segment visited.
This operation is acceptable if applying \(\Bv_0\) to a vector is cheap, and in particular for diagonal $\Bv_0$.
In our case, applying \(\Bv_0\) to a vector costs \(O(n \, \log n)\) operations and quickly becomes time consuming.
Instead, we use an Armijo-like backtracking search, similar to that implemented in TRON \citep[Section~6]{Lin99}.
In the inexact procedure, the Cauchy step must only satisfy the sufficient decrease condition
\begin{equation}\label{eq:cauchyineq}
    m_k(\xC) \leqslant f(\xv_k) + \mu_0\, \gv\T (\xC - \xv_k),
\end{equation}
where $0 < \mu_0 < \tfrac{1}{2}$.
Though the case $\xC = \xv$ satisfies \eqref{eq:cauchyineq}, the backtracking procedure is expected to find a Cauchy point satisfying this criterion before reaching $\xv$.

Our implementation is summarized in \Cref{algo:sclbfgsb}.
In the next sections refer to it as \scaledLBFGSB{}.

\begin{algorithm}\label[algorithm]{algo:sclbfgsb}
    \caption{Overview of the modified L-BFGS-B algorithm.}
    \KwData{$\xv_0$, parameters}
    \For{$k = 0, 1, 2, \dots$}{
    	Identify the binding set $I_+(\xv_k)$ and compute \(\bar \Pv\) using \eqref{eq:pbar}\\
        Compute an inexact Cauchy point $\xv^\mathrm{C}$ along the projected path
        $t \mapsto \mathrm{Proj}\left(\xv_k - t\,\bar\Pv\nabla f(\xv_k)\right)$ \\
        Identify the active face $F_+$\\
        Find a minimizer of the model (\ref{eq:or-qnmodel}) over the affine subspace
        (\ref{eq:activeface}) using preconditioned CG \\
        Perform a strong Wolfe linesearch to find the next iterate $\xv_{k+1}$\\
        Update the L-BFGS operator.
    }
\end{algorithm}

    \subsection{A trust-region Newton method}
    	
Because second-order derivatives are available in our reconstruction problem, we also describe our scaling strategy in the context of a Newton method.

\subsubsection{Presentation of the algorithm}

TRON is a trust-region Newton method proposed by \citet{Lin99}.
As in L-BFGS-B, a general iteration includes the identification of an active face via a Cauchy point and the minimization of a quadratic model over an affine subspace corresponding to the free indices.

The quadratic model (\ref{eq:or-qnmodel}) now uses $\Bv_k = \nabla^2 f(\xv_k)$.
Due to the high cost of evaluating the quadratic model, we only compute an inexact Cauchy point satisfying~(\ref{eq:cauchyineq}).

Computing a Newton direction requires to find an approximate solution of the trust-region subproblem
\begin{equation}\label{eq:tronsubproblem}
    \begin{aligned}
	\min \ & \tfrac{1}{2} \xv\T \Bv_k\xv + \xv\T\gv
	\\ \text{s.t.} \ &
	   \xv \geqslant 0,\quad
	   \forall i \in \mathcal{A}\ x_i = 0, \quad
	   \|\xv -\xv_k\| \leqslant \Delta_k,
	\end{aligned}
\end{equation}
where $\Delta_k > 0$ is the trust-region radius.
To solve this subproblem, the algorithm generates a sequence of minor iterates $\xv^1, \xv^2, \dots$, starting with $\xv^0 = \xC$.
From a minor iterate $\xv^j$, we launch CG on problem \eqref{eq:tronsubproblem}, stopping at convergence or when a point $\xv^j + \wv$ falls out of the feasible set.
In the second case, a projected search (similar to the Cauchy point identification procedure) is used along the path
\begin{equation*}
	t \mapsto \mathrm{Proj}(\xv^j +t\wv)
\end{equation*}
to identify the next minor iterate $\xv^{j+1}$.
Like in the Cauchy procedure, the indices $i$ such that $x^{j+1}_i = 0$ and $w_i < 0$ are added into the active set $\mathcal{A}$, and so the next minor iterate is searched for in a smaller subspace than the current one.
For more information about the subproblem resolution, refer to \citet[Sections 4 and 6]{Lin99}.

The procedure ends when a minor iterate $\xv^j$ satisfies the condition
\begin{equation}\label{eq:minordecrease}
	\|\xv^j - \mathrm{Proj}(\xv^j - \nabla m_k(\xv^j))\|
	\leqslant \varepsilon \|\xv_k - \mathrm{Proj}(\xv_k - \nabla f(\xv_k))\|,
\end{equation}
or when a minor iterate falls outside the trust-region.
Then the decrease of the quadratic model and the decrease of the objective function at the last minor iterate are compared.
Depending on this information, the last minor iterate is accepted as $\xv_{k+1}$ or rejected, and the trust-region radius is updated.

A summary of a TRON iteration is given in \Cref{alg:tron}.

\begin{algorithm}
	\label[algorithm]{alg:tron}
    \caption{Overview of the standard TRON algorithm.}
    \KwData{$\xv_0$, parameters}
    \For{$k = 0, 1, 2, \dots$}{
        Compute an inexact Cauchy Point $\xC$ along the projected path
        $t \mapsto \mathrm{Proj}\left(\xv_k - t\,\nabla f(\xv_k)\right)$\\
        Identify the active set $\mathcal{A}$\\
        $j \leftarrow 0$, $\xv^0 \leftarrow \xC$\\
        \While{ (\ref{eq:minordecrease}) is not satified}
        {
        	Launch CG on problem \eqref{eq:tronsubproblem} to find a direction $\wv$\\
        	Perform a projected search along the path $\mathrm{Proj}(\xv^j + t\wv)$ to obtain $\xv^{j+1}$\\
        	$\mathcal{A} \leftarrow \mathcal{A} \cup \{i\,|\,x^{j+1}_i = 0\mbox{ and } w_i < 0\}$\\
        	$j \leftarrow j+1$\\
        }
        Accept or reject $\xv^j$ as the new iterate and update the trust-region radius\\
    }
\end{algorithm}

\subsubsection{Implementation in Matlab}

In the original TRON Fortran implementation, the conjugate gradient is preconditioned using an incomplete Cholesky factorization of $\Bv_k$.
Such a factorization is not appropriate for large problems because the matrix coefficients are not explicitly available.

We use \bcflash{}, a Matlab implementation of TRON provided by Friedlander and Orban\footnote{Available online at \url{https://github.com/optimizers/NLPLab}}, without preconditioning and where $\Bv_k$ is only used as an operator.

\begin{figure}
    \centering
    \includegraphics[height=.31\columnwidth]{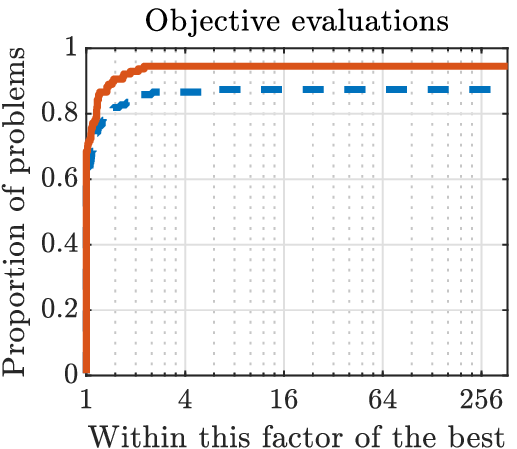}\!
    \includegraphics[height=.31\columnwidth]{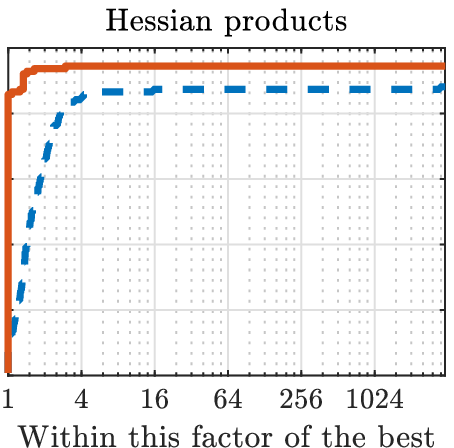}\!
    \includegraphics[height=.31\columnwidth]{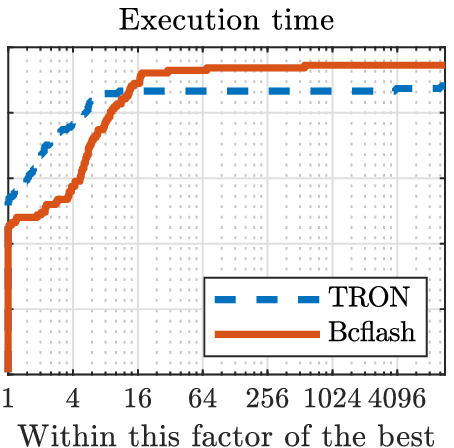}
    \caption{Performance profiles of Bcflash (Matlab) versus TRON without factorization (Fortran).}
    \label[figure]{fig:valTRON}
\end{figure}

For validation, we test \bcflash{} against the Fortran TRON implementation from which the incomplete Cholesky factorization was removed. 
The profiles in \Cref{fig:valTRON} show the performance results on 127 problems from the CUTEst library \citep{Gould15}.
The profiles show that \bcflash{} is more efficient in terms of function evaluations and Hessian products than the Fortran version.
Moreover, the Matlab implementation is more robust.
These results confirm the validity of Bcflash as an implementation of TRON.

\bcflash{} is competitive with the Fortran implementation in terms of execution time, whereas the difference is larger for L-BFGS-B.
This can be partially explained as follows.
In TRON, the bulk of the computation resides in Hessian-vector products.
In both implementations the latter are computed by the CUTEst infrastructure, so this part of the computation is common between them.
In L-BFGS-B, the objective function and gradient are only called at the beginning of the iteration and during the line search.
Moreover, computations related to using and updating the limited-memory operator are difficult to vectorize efficiently, as they include operating on triangular matrices and reordering matrix columns.
Thus, \LBFGSB{} is at a disadvantage because those computations are implemented in Matlab.

A summary of the scaled variant of TRON appears in \Cref{alg:sctron}.
From here on, we refer to it as \scaledbcflash{}.

\begin{algorithm}
	\label[algorithm]{alg:sctron}
    \caption{Overview of the modified TRON algorithm.}
    \KwData{$\xv_0$, parameters}
    \For{$k = 0, 1, 2, \dots$}{
    	Identify the binding set $I_+(\xv_k)$ and compute \(\bar \Pv\) using \eqref{eq:pbar}\\
        Compute an inexact Cauchy Point $\xC$ along the projected path
        $t \mapsto \mathrm{Proj}\left(\xv_k - t\,\bar\Pv\nabla f(\xv_k)\right)$\\
        Identify the active set $\mathcal{A}$\\
        $j \leftarrow 0$, $\xv^0 \leftarrow \xC$\\
        \While{ (\ref{eq:minordecrease}) is not satified}
        {
        	Launch CG preconditioned with $\Pv_{F\!F}$ on problem \eqref{eq:tronsubproblem} to find a direction $\wv$\\
        	Perform a projected search along the path $\mathrm{Proj}(\xv^j + t\wv)$ to obtain $\xv^{j+1}$\\
        	$\mathcal{A} \leftarrow \mathcal{A} \cup \{i\,|\,x^{j+1}_i = 0\mbox{ and } w_i < 0\}$\\
        	$j \leftarrow j+1$\\
        }
        Accept or reject $\xv^j$ as the new iterate and update the trust-region radius\\
    }
\end{algorithm}

%%%%%%%%%%%%%%%%%%%%
\section{Numerical results}\label[section]{S4}
    We now evaluate the performance of \Cref{algo:sclbfgsb} and \Cref{alg:sctron} on two image reconstruction problems in cylindrical coordinates.
For both problems, we compare the performance of \scaledLBFGSB{} and \LBFGSB{}, and that of \scaledbcflash{} and \bcflash{}.
We also compare \scaledbcflash{} with a change-of-variable approach, by using \cflash{}, an implementation of TRON adapted to polyhedral constraints, to solve \eqref{eq:screcmodel}.
In \cflash{}, projections are made onto the feasible set by solving a quadratic problem with Krylov methods \citep{McLaughlin17}.
Performances are compared in terms of projected gradient norm decrease and cumulated number of CG iterations along the reconstruction procedure.

First, to compare the convergence properties of the algorithms, we solve a simplified reconstruction problem, which is quadratic and better scaled than \eqref{eq:recmodel}.
In this first test, we also measure the fraction of execution time spent computing products with \(\Av\) and \(\Cv\) for \bcflash{}, \scaledbcflash{} and \cflash{}, in order to emphasize the high cost of constraints management in the third method compared to that in the two other methods.

In a second test, we use our methods on a real reconstruction problem.
We evaluate the convergence acceleration caused by the use of scaled directions in this more complex case.
To justify the choice of higher-order methods in image reconstruction, we also compare the performance of \scaledbcflash{} with that of a first-order method, the spectral projected gradient (SPG) of \citet{Birgin02}.

In the following tests, we reconstruct images from a $672 \times 1 \times 1160$ synthetic sinogram.
In order to keep reasonable reconstruction times, we only consider 2D images.
The data were created from a XCAT phantom \citep{Segars08} of size $512 \times 512$ in cartesian coordinates.
We reconstruct a discretized image using polar coordinates, with 226 radial subdivisions and 1160 angular subdivisions.
This discretization provides a sufficient resolution to obtain, after conversion, a $512 \times 512$ cartesian image.
Thus, the data creation and the image reconstruction are made using different procedures.
In this problem, $\Av$ has 779,520 rows and 262,160 columns and the initial guess is $\xv_0 = 0$.

Here, choosing another initial guess, such as the filtered back-projection, would bring no improvement in terms of image quality, because we use gradient-based methods for the reconstruction.
During a reconstruction using gradient-based methods (not a Gauss-Seidel method for  instance),  it  is  observed  that  the  low-frequency  components  of  the  image converge  first whereas the  high-frequency  components  (including  edges  and  noise) converge more slowly \citep{Sauer93}.
For this reason, there is no advantage in initializing with filtered back-projection, which is usually noisy.
It could even increase the noise in the final image, whereas beginning with a constant image makes the noise appear only because of the data.
Moreover, the choice of the initial guess it out of the scope of this article, as we focus on convergence speed comparison.

All results below are produced on an Intel\textsuperscript{\tiny\textregistered} Xeon\textsuperscript{\tiny\textregistered} E5-2637 v4 processor at 3.50~GHz and 32~GB of RAM.

\subsection{Simplified problem}

We first consider the regularized linear least-squares problem
\begin{equation}\label{eq:testpb1}
	\min_{\xv \geqslant 0} \tfrac{1}{2} \|\Av\xv - \bv\|^2 + \tfrac{1}{2}\lambda \|\Kv \xv\|^2,
\end{equation}
where $\Kv\xv$ models the differences between adjacent pixels of $\xv$.
In this simplified reconstruction problem, we drop the weight matrix $\Vv$, which is equivalent to assuming that all attenuation measures have the same variance,
and we choose a $L_2$ regularization to keep the problem quadratic. 
We set $\lambda = 10^{-2}$ because it provides reasonable image quality and convergence speed.

\begin{figure} %[ht]
    \centering
    \includegraphics[height=.39\columnwidth]{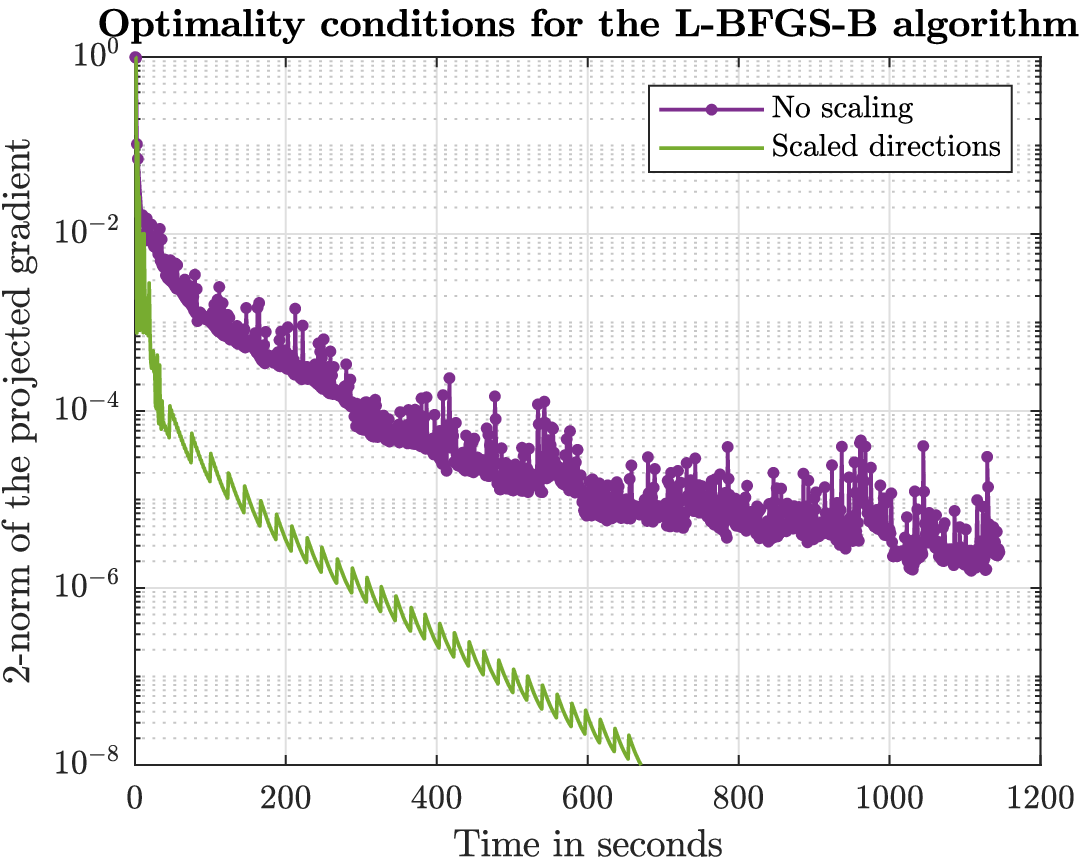}
    \includegraphics[height=.39\columnwidth]{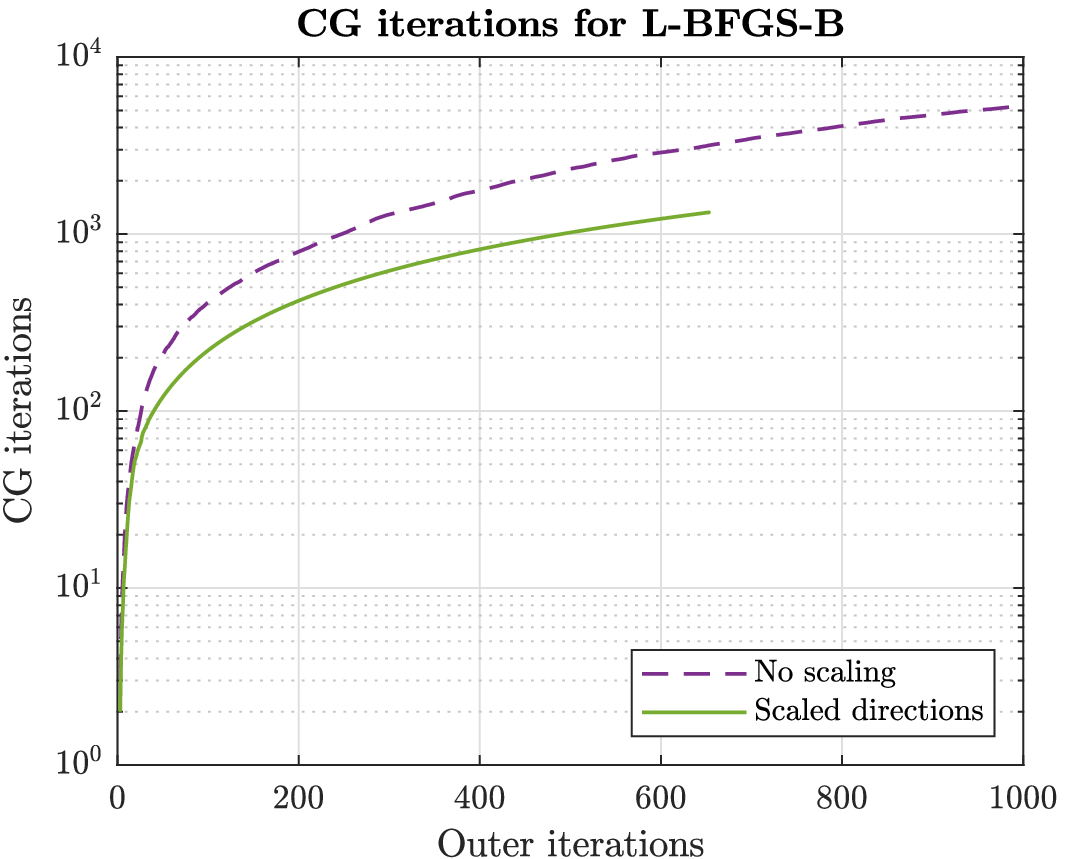}
    \caption{Convergence results for L-BFGS-B on (\ref{eq:testpb1}).}
    \label[figure]{fig:lbfgsb1}
\end{figure}

\Cref{fig:lbfgsb1} shows the comparison between \LBFGSB{} and \scaledLBFGSB{}.
The left and right plots compare the decrease of the optimality residual and the cumulated number of CG iterations, respectively.
Like \citet[Section~5.2]{Byrd95}, we set the CG relative stopping criterion to $\min(0.1,  \sqrt{\|r^c\|})$, where $r^c$ is the residual at the beginning of the CG procedure.
\citet{Conn88} recommend this choice to enforce superlinear convergence.
We observe that the projected gradient norm decreases faster in the scaled case, especially in the first iterations, and that \scaledLBFGSB{} requires about half as many CG iterations per outer iteration as \LBFGSB{}.
The use of a nondiagonal \(\Bv_0\) yields L-BFGS approximations that are closer to the problem Hessian and lead to better progress at each step.
For this reason, we see on the left plot that \scaledLBFGSB{} decreases the projected gradient norm more than \LBFGSB{} while doing less outer iterations.
Moreover, each outer iteration has a lower cost in terms of CG iterations due to the use of a preconditioner when solving~\eqref{eq:quadinF}.

However, the performance of \scaledLBFGSB{} is not sufficient.
Though much progress is achieved during the 30 first seconds, the convergence seems to switch to a linear behavior at some point, and it takes more than two minutes to decrease the projected gradient norm by a factor of~$10^5$.

\begin{figure} %[ht]
    \centering
    \includegraphics[height=.38\columnwidth]{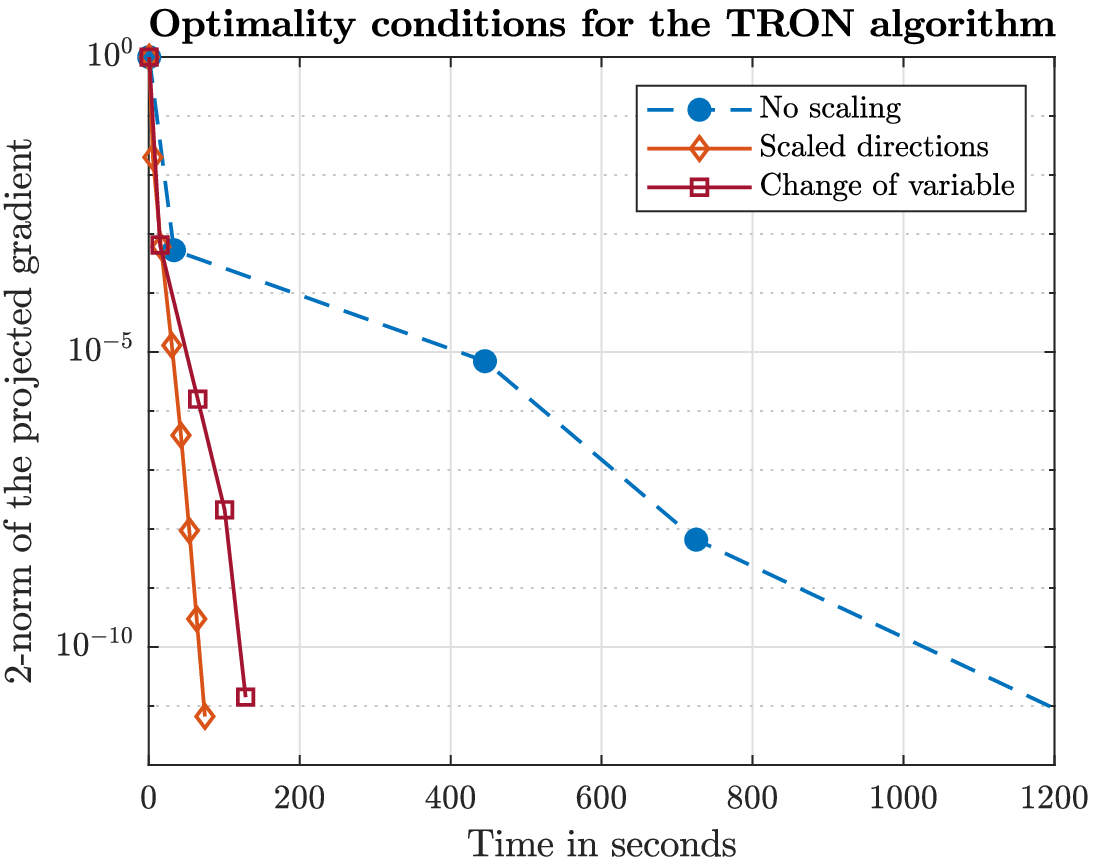}
    \includegraphics[height=.38\columnwidth]{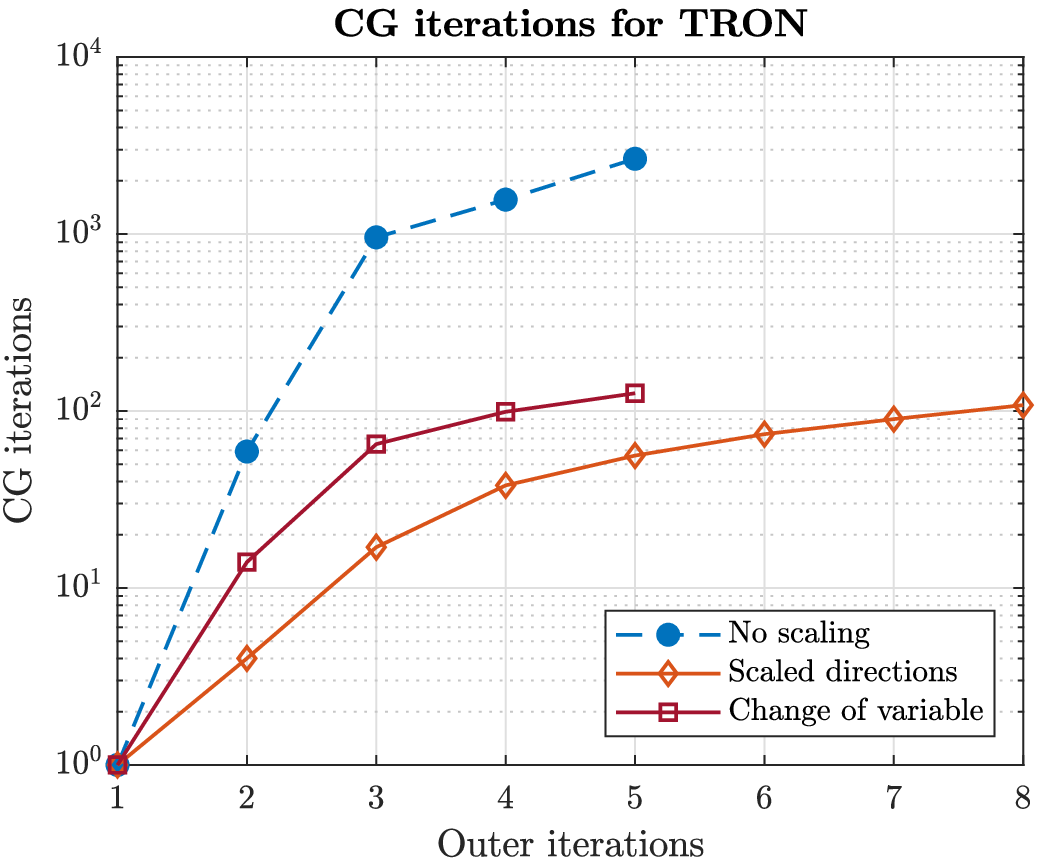}
    \caption{Convergence results for TRON on (\ref{eq:testpb1}). Each marker represents an outer iteration of the Newton method.}
    \label[figure]{fig:tron1}
\end{figure}

\Cref{fig:tron1} reports corresponding results for TRON, where we compare \bcflash{}, \cflash{} and \scaledbcflash{}.
In this test, we set the CG stopping relative tolerance to $10^{-3}$, which gave the best results for the two scaled methods.
On this figure, each marker stands for one outer iteration of the Newton method.
Both \cflash{} and \scaledbcflash{} decrease the projected gradient much faster than \bcflash{}.
Even though \scaledbcflash{} requires more iterations than \cflash{} for the same gradient decrease, it converges faster.

\Cref{tab:stat1} gathers some statistics about the execution of TRON.
Due to the problem size, a significant part of the execution time is associated with products with $\Av$ or its transpose.
We see in the table that this time fraction is similar for \bcflash{} and \scaledbcflash{}, both of which work in the original space, while it is much smaller for \cflash{}, which works in the scaled space.
This difference is associated with the time spent computing orthogonal projections onto affine constraints in \cflash{}.
Indeed, $23\%$ of the solve time is spent computing products with $\Cv$ or its transpose, most of which are computed during orthogonal projections.
Note that the time spent computing product with $\Av$ or its transpose are approximatey identical for \scaledbcflash{} and \cflash{}.
In \scaledbcflash{}, time is saved on projections while the cost of conjugate gradient iterations remains similar.

\begin{table} %[ht]
    \centering
    \begin{tabular}{|l|ccc|}
        \hline
        TRON variant & \bcflash{} & \scaledbcflash{} & \cflash{} \\
        \hline
        Total execution time & 1235~s & 76~s & 131~s \\
        Time spent doing $\Av$-products & 1203~s & 73~s & 73~s \\
        Time spent doing $\Cv$-products & 0~s & 0.7~s & 30~s \\
        Time fraction for $\Av$-products & 97~\% & 95~\% & 56~\% \\
        Time fraction for $\Cv$-products & 0~\% & 1~\% & 23~\% \\
        \hline
        Outer iterations & 4 & 7 & 4 \\
        Avg. time per outer iteration & 308~s & 11~s & 32~s \\
        CG iterations & 2663 & 108 & 126 \\
        Avg. time per CG iteration & 0.46~s & 0.7~s & 1.0~s \\
        \hline
    \end{tabular}
    \caption{Execution statistics for the three versions of the TRON algorithm: fraction of time spent doing products with $\Av$ and $\Cv$}
    \label[table]{tab:stat1}
\end{table}

\subsection{Reconstruction problem}

In the second test, we compare \Cref{algo:sclbfgsb} and \Cref{alg:sctron} on the reconstruction problem
\begin{equation}\label{eq:testpb2}
	\min_{\xv \geqslant 0} \tfrac{1}{2} \|\Av\xv - \bv\|^2_\Vv + \lambda \phi(\Kv \xv),
\end{equation}
where 
$\phi: \qv \mapsto \sum_{i} \sqrt{\delta^2 + q_i^2}$ is an edge-preserving $L_2 / L_1$ penalty with $\delta > 0$,
and $\Vv$ is the statistical weight matrix defined in \eqref{eq:recmodel}.

Problem (\ref{eq:testpb2}) should be more difficult to solve than (\ref{eq:testpb1}), even for scaled methods.
Indeed, the addition of weights deteriorates the Hessian conditioning, and the penalty is not quadratic.
The objective Hessian in (\ref{eq:testpb2}) has the form
\begin{equation}\label{eq:hessianpb2}
	\nabla^2 f(\xv) = \Av\T \Vv \Av + \Kv\T\Nv(\Kv\xv)\Kv,
\end{equation}
where $\Nv(\Kv\xv)$ is a diagonal matrix with general term
\begin{equation*}
    n_{ii} = \delta^2 / \left( \delta^2 + [\Kv \xv]_i^2 \right)^{3/2}.
\end{equation*}
This Hessian is not block-circulant due to the presence of $\Vv$ and $\Nv(\kv\xv)$.
To compute $\Pv$, we block-diagonalize a block-circulant approximation of \eqref{eq:hessianpb2}.
First, we define $\hat{\Vv}$ by averaging the diagonal blocks of $\Vv$.
If $\Vv = \mathrm{diag}(\Dv_1, \cdots, \Dv_{n_b})$, where $\Dv_1, \cdots, \Dv_{n_b}$ are diagonal blocks, we create 
\begin{equation*}
	\hat{\Vv} = \mathrm{diag}(\hat\Dv, \cdots, \hat\Dv),
	\quad\mbox{where}\quad
	\hat\Dv = (\Dv_1 + \cdots + \Dv_{n_b}) / n_b.
\end{equation*}
In addition, instead of $\Nv(\Kv\xv)$, we use $\Nv(0) = (1 / \delta) \Iv$.
Finally the scaling matrix is based on the block-diagonalization of $\Av\T \hat\Vv \Av + \Kv\T \Nv(0)\Kv$.
Due to this intermediate approximation we expect the scaling of problem~\eqref{eq:testpb2} to be less efficient than that of~\eqref{eq:testpb1}.

We choose $\lambda = 10^{-4}$ and $\delta = 10^{-1}$.
These parameters results from a trade-off between the speed of convergence and the quality of reconstructed images, namely the noise level and the attenuation of sharp edges compared to the original image \cite[section 4.5]{Hamelin09d}.

\begin{figure}
    \centering
    \includegraphics[height=.39\columnwidth]{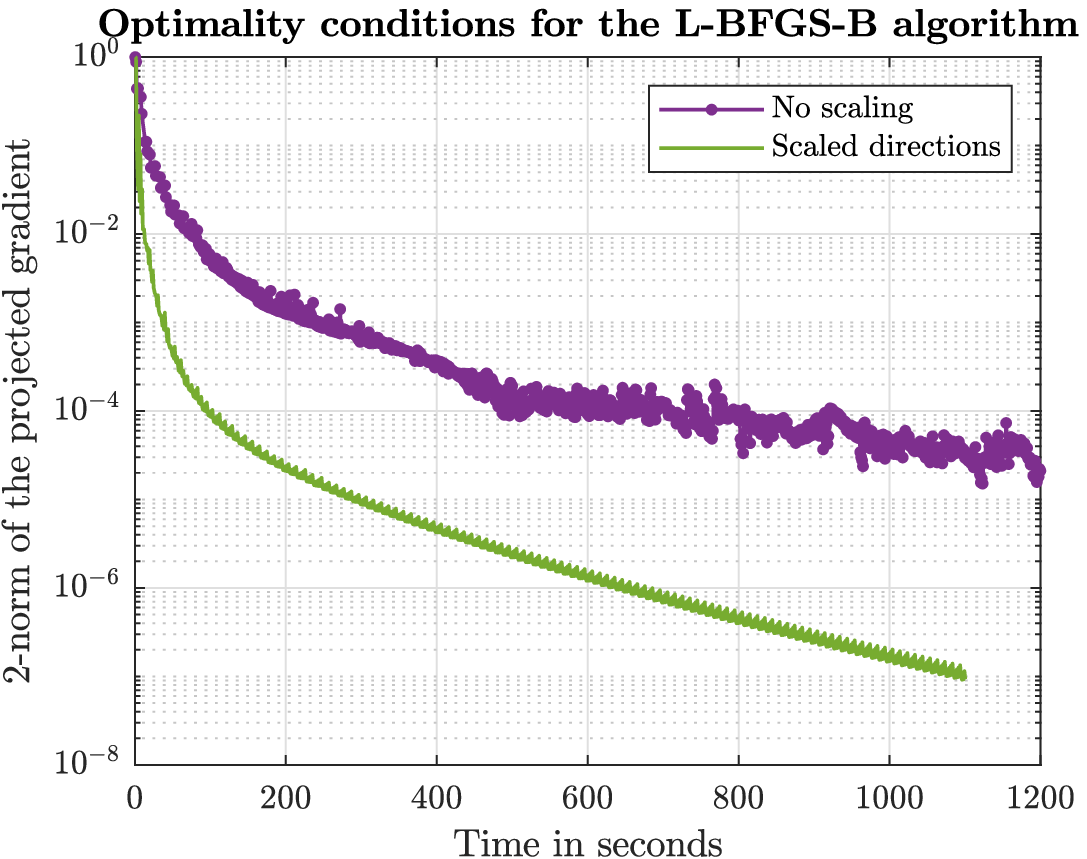}
    \includegraphics[height=.39\columnwidth]{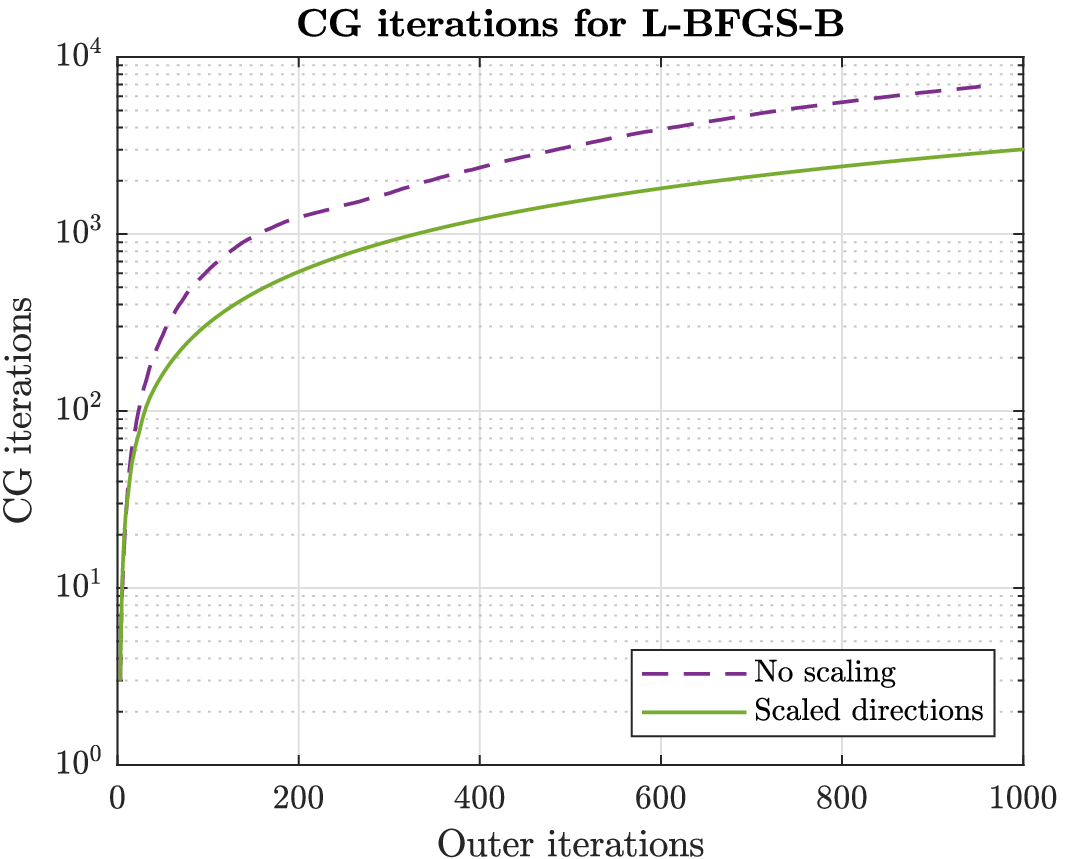}
    \caption{Convergence results for L-BFGS-B on (\ref{eq:testpb2})}
    \label[figure]{fig:lbfgsb2}
\end{figure}

\begin{figure}
    \centering
    \includegraphics[height=.39\columnwidth]{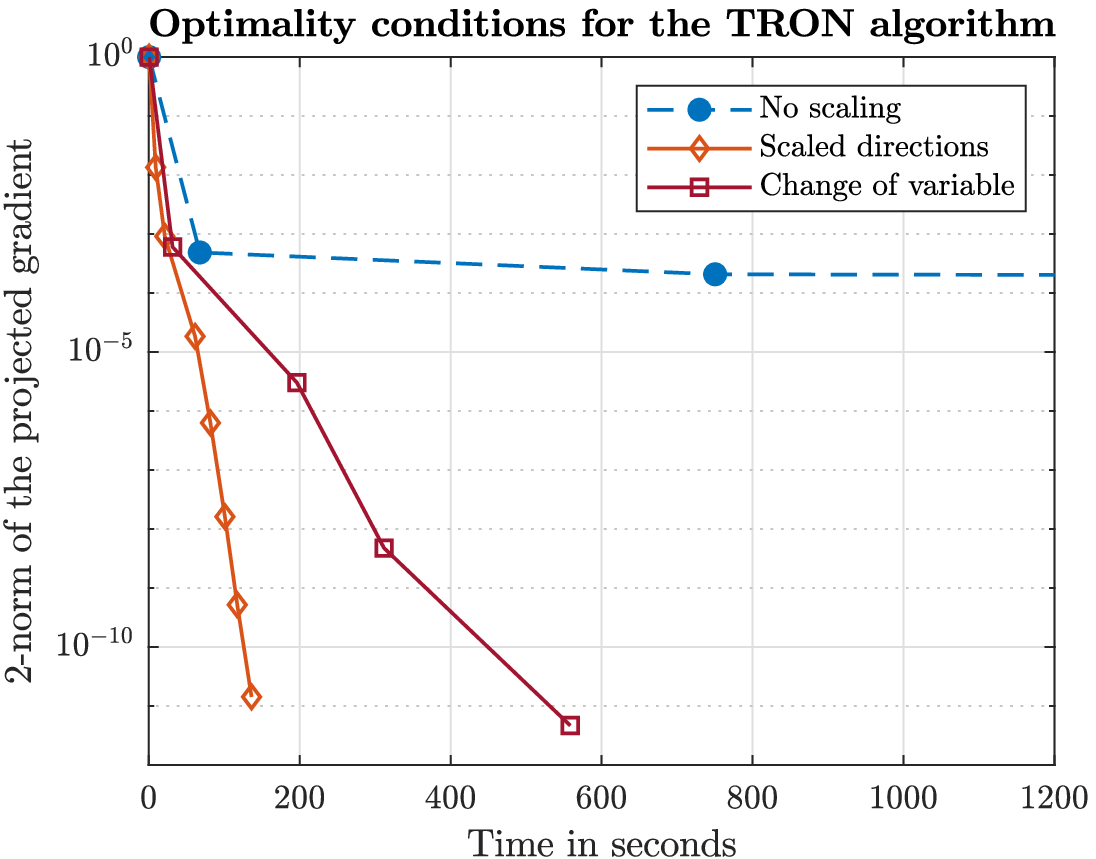}
    \includegraphics[height=.39\columnwidth]{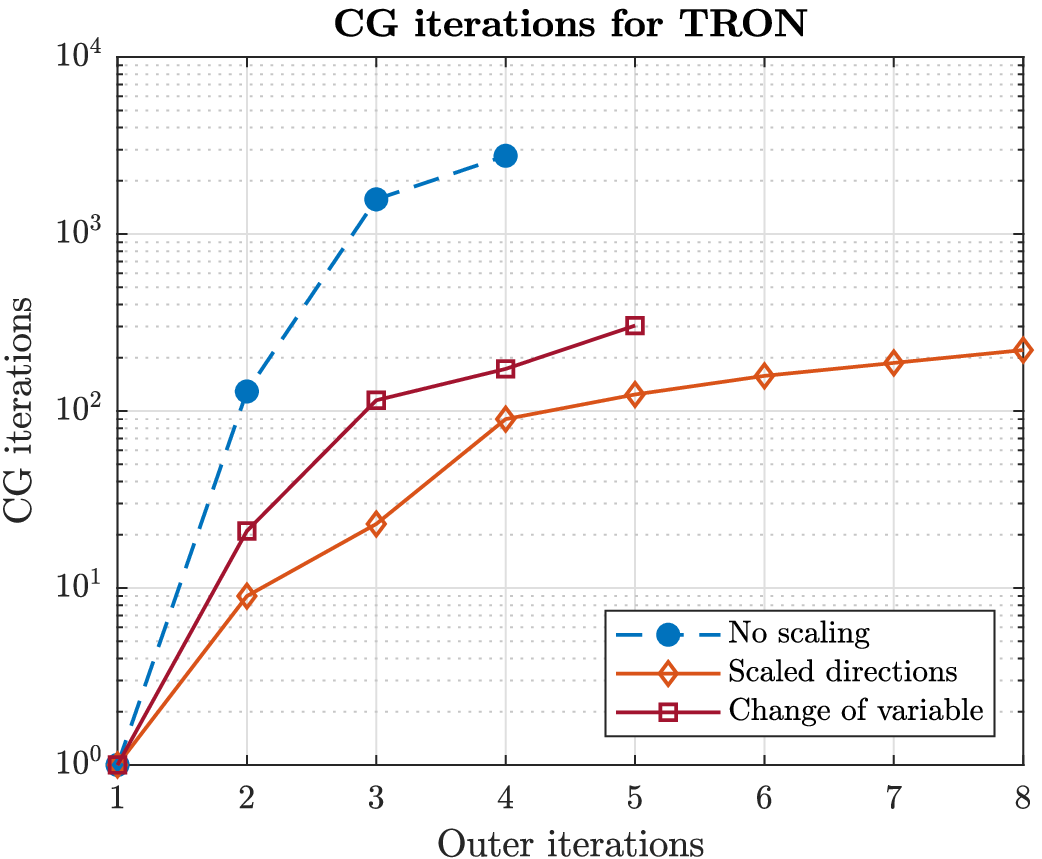}
    \caption{Convergence results for TRON on (\ref{eq:testpb2}). Each marker represents an outer iteration of the Newton method.}
    \label[figure]{fig:tron2}
\end{figure}

\Cref{fig:lbfgsb2} and \Cref{fig:tron2} show convergence results for L-BFGS-B and TRON, respectively.
The solve time is longer than on~(\ref{eq:testpb1}) for all solvers, and the impact of the scaled solver is not as pronounced as in~(\ref{eq:testpb1}).
The \scaledLBFGSB{} decreases the projected gradient by a factor of $10^4$ about $9$ times faster than \LBFGSB{} on~(\ref{eq:testpb1}), but only $5$ times faster on~(\ref{eq:testpb2}).

In the case of TRON, the advantage of \scaledbcflash{} over \cflash{} is larger on \eqref{eq:testpb2} than on \eqref{eq:testpb1}.
The \scaledbcflash{} decrease the projected gradient by a factor of $10^7$ about $1.8$ times faster than \cflash{} on~(\ref{eq:testpb1}), and $3.1$ times faster on~(\ref{eq:testpb2}).
So \scaledbcflash{} is less affected by the scaling deterioration than \cflash{}.
Note that the unscaled \bcflash{} does not manage to achieve a better gradient decrease than $10^{-4}$.
This is due to the strict CG tolerance, $10^{-3}$, that we use for the three solvers.
With a larger tolerance, such as $10^{-1.5}$, \bcflash{} solves the problem, but after a much longer time than te two scaled versions.

\begin{figure}
    \centering
    \includegraphics[height=.39\columnwidth]{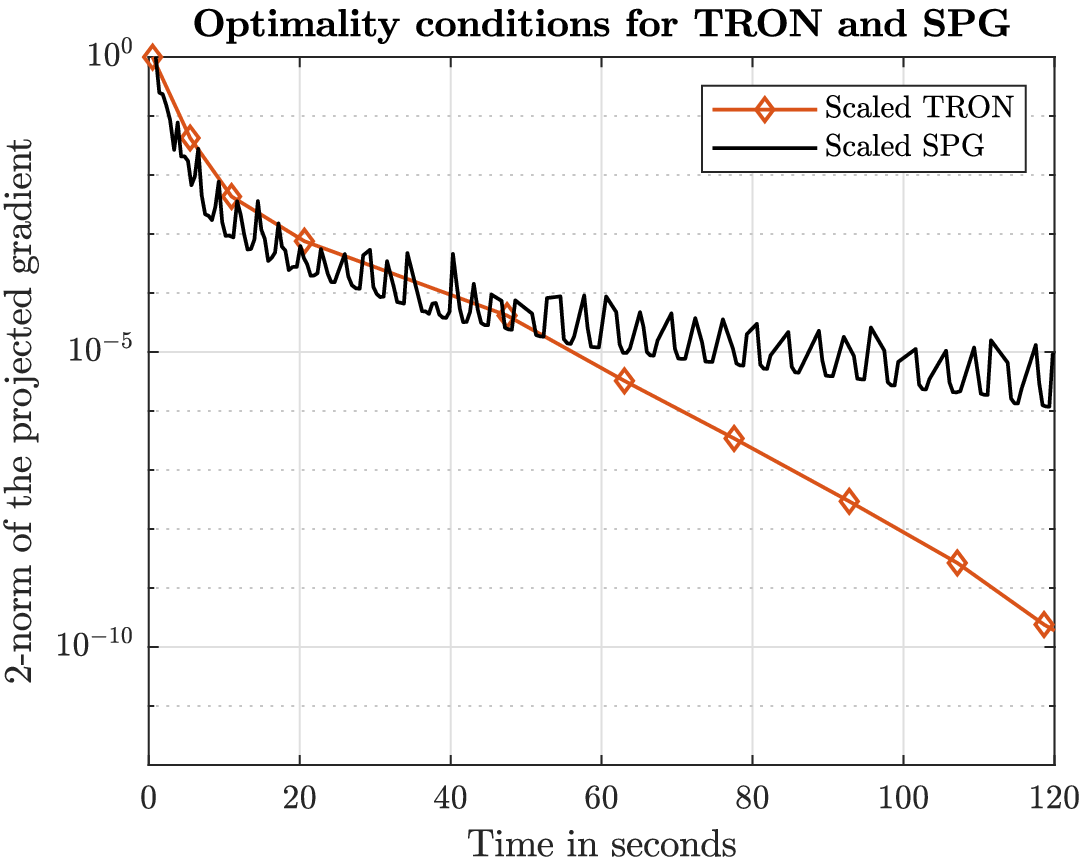}
    \includegraphics[height=.39\columnwidth]{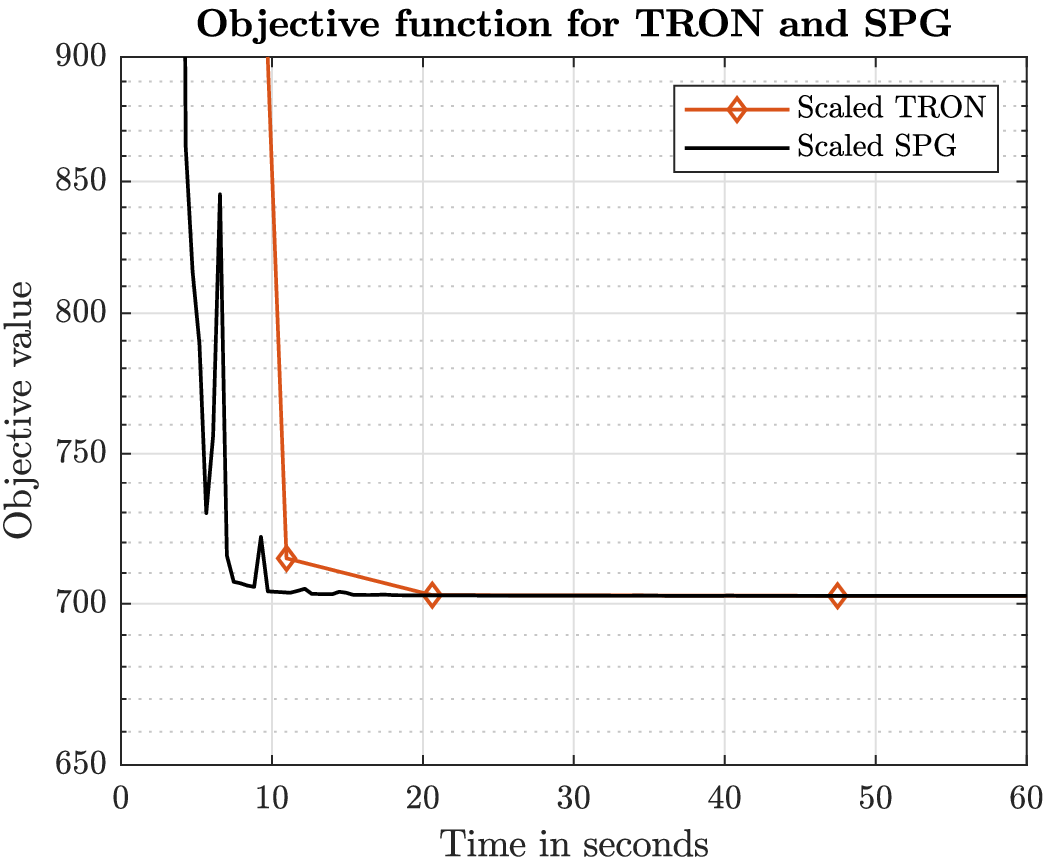}
    \caption{Comparison of TRON and SPG on (\ref{eq:testpb2})}
    \label[figure]{fig:tronspg}
\end{figure}

\Cref{fig:tronspg} shows a comparison between the convergence of \scaledbcflash{} and the spectral projected gradient (SPG) of \citet{Birgin02}, which is appealing because the cost of each iteration is low.
We use a Matlab implementation of SPG\footnote{Available online at \url{https://github.com/optimizers/NLPLab}},
which we modify to employ the same scaling strategy as \Cref{algo:sclbfgsb} and \Cref{alg:sctron}.
The CG tolerance for \scaledbcflash{} is $10^{-2}$.
The resulting algorithm is close to that of \citet{Bonettini08}, except that the scaling is nondiagonal in our case.
SPG decreases the objective function faster than \scaledbcflash{} in the first iterations, which is to be expected from a first-order method.
However, after a few iterations, the SPG projected gradient norm decreases slowly, whereas the decrease rate is superlinear for TRON.
With a 2 minute time limit, \scaledbcflash{} decreases the projected gradient by a factor of $10^{9}$ while SPG only achieves a reduction of about $10^{6}$.
We conclude from the results that \scaledbcflash{} is a more appropriate for solving the reconstruction problem at tolerances stricter than $10^{-5}$.

These numerical results show that the scaling strategy brings the expected performance improvements for the problems we are interested in.
In the case of TRON, this improvement is better than that we obtain with a change of variable, due to the high cost of orthogonal projections in the second approach.

%%%%%%%%%%%%%%%%%%%%
\section{Conclusion}

We presented a scaling strategy for bound-constrained problems inspired from a change of variable, and integrated it into two projected-directions algorithms, L-BFGS-B and TRON.
Though, our strategy can be implemented into most projected-directions algorithms for large bound-constrained problems with little code modifications.
In this paper, we adopted a practical point of view, as we gave details about the implementation of scaling for each subroutine of the algorithms, including the preconditioning of CG to solve quadratic subproblems that appear in higher-order methods.
The numerical tests on badly scaled image reconstruction problems show that this approach gives better results than a change of variable, especially because the management of contraints is cheaper.

These results are promising for applications in X-Ray CT reconstruction, as they show the feasibility of reconstructing images in cylindrical coordinates.
The partially diagonal scaling ensures the efficiency of the procedure, making fast and memory-saving reconstructions possible.
In particular, TRON is a good candidate for applications which require to solve the reconstruction problem with a tight tolerance.
As TRON can solve \eqref{eq:testpb2} with tolerance \(10^{-10}\), the reconstructed image is very close to the problem solution, and can be used as a reference to evaluate the convergence speed of other algorithms.

Here we have implemented the scaling into generic methods for large bound-constrained problems and solved generic reconstruction problems.
In a future work, we will produce scaled methods for specific applications in medical image reconstruction, in order to combine the memory savings provided by the cylindrical coordinates with performance compatible with clinical applications.
In particular, cylindrical coordinates are appropriate when the source and the detector follow a circular trajectory around the object to investigate, like in cone-beam computed tomography or in nondestructive testing.
Structured system matrices can appear for other acquisition protocols, as long as the image discretisation yields geometrical invariances.
Thus, block-circulant system matrices may also appear in helical computed tomography by using an helical discretization.

Circulant structures also appear in other imaging problems, for which our strategy can be applied.
In general, our methods can prove to be useful in applications which lead naturally to non-diagonal scaling operators, like partial differential equations and optimal control, or when the optimization problem is too badly scaled for a diagonal scaling to be efficient.

%%%%%%%%%%%%%%%%%%%%%%%%%%%%%
% BOTTOM
%%%%%%%%%%%%%%%%%%%%%%%%%%%%%

\small
\bibliographystyle{abbrvnat}
%\bibliography{biblio/revueabr,biblio/mois_ct_en,biblio/baseLION}   

\normalsize

\end{document}